\documentstyle{czjphys209} 
 
\RequirePackage{amssymb}
\newtheorem{Theorem}{Theorem}[section]
\begin{document}

\title{Quantized Dirac Operators}

\authori{Hans Plesner Jakobsen}

\addressi{Mathematics Department, 
University of Copenhagen, Universitetsparken 5, 
DK-2100 Copenhagen \O, Denmark}

\authorii{}

\addressii{}

\authoriii{}

\addressiii{}

\headtitle{Quantized Dirac\ldots}

\headauthor{Jakobsen} 

\specialhead{Jakobsen: The Quantized Dirac Operator}


\evidence{A}

\daterec{XXX}    

\cislo{0}  \year{2000}

\setcounter{page}{1}

\pagesfromto{000--000}



\maketitle

\begin{abstract}
We determine what should correspond to the Dirac operator on certain
quantized hermitian symmetric spaces and what its properties are. A new insight
into the quantized wave operator is obtained. 
\end{abstract}

\section {Introduction}

In the articles \cite{j1},  \cite{hj1}, \cite{hj2},  \cite{j2},
\cite{j3}, and \cite{j4} it was shown how intertwining differential
operators on spaces of vector valued holomorphic functions on
hermitian symmetric spaces of the non-compact type are intimately
connected both with homomorphisms between generalized Verma modules
and with singular unitary representations through what was called the
``missing $k$ types''. It was clear from the outset that the
construction was applicable in a much broader context for classical
(semi-simple) Lie algebras. Later, it was natural to extend it further
to quantum groups, in particular after the classification of
unitarizable highest weight modules was obtained \cite{j5}. 

In connection with the generalization to quantum groups of the above,
one needs of course to find good replacements of hermitian symmetric
spaces, holomorphic functions, and, finally, differential operators.
We believe that we through the discovery of quadratic algebras
associated to hermitian symmetric spaces, see e.g. \cite{j6}, have
found the right setting. 

In a series of papers, c.f. \cite{d1}, \cite{d2}, and references cited
therein, Dobrev has studied a quantization of the d'Alembert operator
and its invariance properties and has also addressed the issue of
relevant spaces. The quantized wave operator is also discussed in
\cite{j7} which has some overlap with Dobrev's work, both mathematical
and notational, but in fact, Dobrev's methods have much more in common
with the earlier work cited above. Even more so, the methods in those
references reach much further as we will illustrate below by giving
``the quantization of the Dirac operator'' in the setting of quantum
groups. Our approach shows that the quantized wave operator is a major
part of the quantized Dirac operator; in some sense the term that
represents the most ``non-commutative'' part.

\section{Basics}

For any hermitian symmetric space there is a (matrix valued) first
order holomorphic differential operator that intertwines two unitary
highest weight representations. Sometimes there are even two. These
operators might be called ``Dirac operators'' and they all have
quantizations. In order to keep the presentation simple and at the
same time treat the most relevant case for physics, we will only give
the details here in the case of $su(2,2)$ and its quantized enveloping
algebra ${\mathcal U}_q(su(2,2))$. To be specific, ${\mathcal U}_q={\mathcal U}_q(su(2,2))$ is
generated by 
\begin{equation}
E_\mu, E_\nu, E_\beta, F_\mu, F_\nu, F_\beta, K^{\pm1}_\nu,
K^{\pm1}_\mu,  K^{\pm1}_\beta,
\end{equation} 
with the usual relations of which we only list those we need in the
following:

\begin{eqnarray}\nonumber
  F_\nu^2F_\beta +(q+q^{-1})F_\nu F_\beta F_\nu +F_\beta F^2_\nu=0,
&& \kern-1em F_\mu^2F_\beta +(q+q^{-1})F_\mu F_\beta F_\nu +F_\beta
F^2_\mu=0,\\\nonumber
F_\beta^2F_\nu +(q+q^{-1})F_\beta F_\nu F_\beta +F_\nu F^2_\beta=0,
&& \kern-1em F_\beta^2F_\mu +(q+q^{-1})F_\beta F_\mu F_\beta +F_\mu
F^2_\beta=0,\\
F_\nu F_\mu= F_\mu F_\nu,\; K_\nu F_\beta=q F_\beta K_\nu,&& \kern-1em K_\mu
F_\beta =q F_\beta K_\mu,\; K_\beta F_\beta=q^{-2} F_\beta
K_\beta,\\\nonumber 
 \left[E_\beta,F_\beta\right]=\frac{K_\beta-K_\beta^{-1}}{q-q^{-1}}, &&
 \left[E_\mu,F_\mu\right]=\frac{K_\mu-K_\mu^{-1}}{q-q^{-1}}.
\end{eqnarray}

\begin{eqnarray}
\bigtriangleup (E_i) &=&E_i\otimes 1+K_i\otimes E_i,\,\bigtriangleup
(F_i)=F_i\otimes K_i^{-1}+1\otimes F_i,  \label{coprod1} \\
\bigtriangleup (K_i) &=&K_i\otimes
K_i,\,S(E_i)=-K_i^{-1}E_i,\,S(F_i)=-F_iK_i,\,S(K_i)=K_i^{-1},\textrm{ and } 
\nonumber \\
\epsilon (E_i) &=&\epsilon (F_i)=0,\,\epsilon (K_i)=1.  \nonumber
\end{eqnarray}

\subsection{The quadratic algebra}

Let $w_1=F_\beta, w_2=F_\mu F_\beta- qF_\beta F_\mu, w_3=F_\nu
F_\beta- q F_\beta F_\nu, w_4=F_\mu w_3- qw_3 F_\mu$. Then the four
$w_i$s generate a quadratic algebra denoted ${\mathcal A}_q$. The
relations are the usual ones,
\begin{eqnarray}
  w_1w_2=qw_2w_1,  w_1w_3=qw_3w_1,  w_3w_4=qw_4w_3,\\\nonumber
  w_2w_4=qw_4w_2, 
  w_2w_3=w_3w_3, w_1w_4-w_4w_1=(q-q^{-1})w_2w_3.
\end{eqnarray}

The element $w_1w_4-qw_2w_3$ is in the center of ${\mathcal A}_q$, and
at a generic $q$, the center is generated by this element.

There is an action (this is the co-adjoint action $k\star w= \sum_ib_i
wS(a_i)$ if $\triangle k= \sum_i a_i\otimes b_i$) of the algebra ${\mathcal U}_q(su2\times su(2))$ on
  ${\mathcal A}_q$ given by 
\begin{eqnarray}\nonumber
    F_\mu \star w_1= w_2,\;F_\mu \star w_3= w_4,&&F_\mu \star
    w_2=0,\;F_\mu \star w_4=0,\;\\
E_\mu \star w_2= w_1,\;E_\mu \star w_4= w_3,&&E_\mu \star
    w_1=0,\;F_\mu \star w_3=0,\;\\\nonumber
 K_\mu \star w_1= qw_1,\;K_\mu \star w_3= qw_3,&&K_\mu \star
    w_2=q^{-1}w_2,\;K_\mu \star    w_4=q^{-1}w_4.
  \end{eqnarray}

This action is derived from a left action in the quantized enveloping
algebra but is not identical with this action. For later use we
record here the identities $F_\mu w_2=q^{-1}w_2 F_\mu$ and $E_\mu
w_2=w_1 K^{-1}_\mu$.

There are analogous equations for the index $\nu$ obtainable simply by
the interchange $2\rlarrows3$.

Observe that 

\begin{equation}
  w_4^Nw_1=w_1w_4^N-(q-q^{-2N+1})w_2w_3w_4^{N-1}.
\end{equation}

Later on, we shall use functions defined on the quadratic algebra
${\mathcal A}_q$, in particular polynomial functions. Actually, we
wish to work with functions on ${\bm C}^4$. For this reason we
introduce the symbols 
\begin{equation}
  [{\bm\gamma}]_q!=
  [\gamma_1]_q![\gamma_2]_q![\gamma_3]_q![\gamma_4]_q!,  
\end{equation}
where fore any non-negative integer $n$, $[n]_q!$ denotes the usual
$q$ factorial.

Likewise,
$w^{{\bm\gamma}}=w_1^{\gamma_1}w_2^{\gamma_2}w_3^{\gamma_3}w_4^{\gamma_4}$
in ${\mathcal A}_q$, and if $z=(z_1,z_2,z_3,z_4)\in {\bm C}^4$ then 
$z^{{\bm\gamma}}=z_1^{\gamma_1}z_2^{\gamma_2}z_3^{\gamma_3}z_4^{\gamma_4}$.

To a (polynomial)  function $f$ on ${\mathcal A}_q$ we now introduce
the function $\Psi_f$ on ${\bm C}^4$ by
\begin{equation}\label{trans}
  \Psi_f(z)=f\left(\sum_{{\bm\gamma}}\frac{z^{{\bm\gamma}}}{ [{\bm\gamma}]_q}w^{{\bm\gamma}}\right).
\end{equation}

The motivation for this normalization is perhaps a little flimsy, but
we mention that that the so-called ``divided powers'' algebra is more
fundamental than others. Furthermore, the space of functions thus
defined on ${\bm C}^4$ actually comes equipped with an associative
$\star$-product. See \cite{j7} for further details.

Later on we shall encounter right multiplication operators on the
space of functions on ${\mathcal A}_q$, i.e. operations of the form
$(({w_0^R})^\dagger\cdot f)(w^{{\bm\gamma}})=f(w^{{\bm\gamma}}w_0)$.
Below, $w_0$ will be either $w_1,w_2,w_3,w_4$, or $w_1w_4-qw_2w_3$.
Transformed into operators on functions on ${\bm C}^4$ via
(\ref{trans}) these become $q$-differential operators expressible in
terms of, among other things, the usual $q$-differential operators on
${\bm C}^1$ as well as  scaling operators (c.f. \cite{j7}, \cite{d1}),
\begin{eqnarray}
  (\left[{\frac{\partial}{\partial
    x}}\right]_q\psi)(x)&=&\frac{f(xq)-f(xq^{-1})}{q-q^{-1}},\textrm{ and}\\
\label{44}K_i(z_1^{\alpha_1}z_2^{\alpha_2}z_3^{\alpha_3}z_4^{\alpha_4})&=&q^{-\alpha_i}z_1^{\alpha_1}z_2^{\alpha_2}z_3^{\alpha_3}z_4^{\alpha_4}\quad
  i=1,2,3,4.
\end{eqnarray}

In terms of un-quantized operators, if $q=e^{\hbar}$ then
$K_i=\sum_n (-\hbar S_i)^n$,  where $S_i=z_i\frac{\partial}{\partial z_i}$.
Also observe that $K_i=\left[\frac{\partial}{\partial z_i}\right]_qz_i-q
  z_i\left[\frac{\partial}{\partial z_i}\right]_q$.

We obtain:
\begin{eqnarray}
  (w_4^R)^\dagger=\left[\frac{\partial}{\partial z_4}\right]_q,
  \;(w_2^R)^\dagger=K_4\left[\frac{\partial}{\partial z_2}\right]_q,
  \;(w_3^R)^\dagger=K_4\left[\frac{\partial}{\partial z_3}\right]_q,
  \textrm{ and } \\
(w_1^R)^\dagger=K_2K_3K_4^2\left[\frac{\partial}{\partial
  z_1}\right]_q+z_4(1-q^{-2})K_4 \square_q
\textrm{, where } \\\square_q \stackrel{Def.}
{=}(w_1w_4-qw_2w_3)^\dagger=K_2K_3\left[\frac{\partial}{\partial
    z_1}\right]_q\left[\frac{\partial}{\partial
    z_4}\right]_q-q\left[\frac{\partial}{\partial
    z_2}\right]_q\left[\frac{\partial}{\partial z_3}\right]_{q}.
\end{eqnarray}

The last defined operator is of course nothing else but the quantized
wave operator. At the most singular point of unitarity for
representations in spaces of scalar valued functions it is an
intertwining differential operator.

\section{Induced representations}

Let ${\mathcal B}^+$ denote the part of ${\mathcal U}_q$ generated by
the elements $K_i^{\pm1}, E_\beta,E_\nu, E_\mu$, and let ${\mathcal
  B}^-$ be defined analogously. Let $\chi_\lambda$ be a 1-dimensional
representation of ${\mathcal B}^+$ corresponding to the weight
$\lambda$. Let
\begin{equation}
H^0(\chi_\lambda)=\{f:{\mathcal U}_q\mapsto {\bm C}\mid
f(ub^+)=\chi^{-1}_\lambda(b^+)f(u) \textrm{ and }f\in \textrm{FIN}({\mathcal
    U}_q(su(2)\times su(2)))\}
\end{equation}
with the action $(x\cdot f)(u)= f(S(x)u)$, $S$ being the antipode. The
condition FIN denotes the usual finiteness condition.

Clearly, any element $u_0\in {\mathcal U}_q$ may in some sense be
viewed as an intertwiner through the map $f\mapsto f_{u_0}$ where 
\begin{equation}
  \label{eq:la}
  f_{u_0}(u)=f(u\cdot u_0),
\end{equation}
but if we further want the target space
to be of the form $H^0(\chi_{\lambda_1})$ for some other 1-dimensional 
representation of ${\mathcal B}^+$, we can assume
that $u_0\in {\mathcal B}^-$ (in fact ${\mathcal N}^-$) and we must
then clearly have that
\begin{equation}
 \forall \;b^+,u:f(u\cdot b^+\cdot
 u_0)=\chi^{-1}_{\lambda_1}(b^+)f(u\cdot u_0).
\end{equation}
This, on the other hand can be directly translated into the following
statement:

Let $M_\lambda$ be the Verma module of highest weight $\lambda$ and
let $0\neq v_\lambda$ denote the highest weight vector. Then $u_0\cdot
v_\lambda$ must define a primitive weight vector in $M_\lambda$.

\begin{Theorem}[H-J, \cite{hj1}, \cite{hj2}, \cite{j4} ]
  \label{eq:t1}
  There is a bijective correspondence between intertwining
  differential operators and homomorphisms between generalized Verma
  modules. 
\end{Theorem}

\section{Dirac operators}

In the present case we shall study two representations and
intertwiners (which will be the quantized Dirac operators) thereon.

Suppose $\chi^1_\lambda(K_\mu)=q$, $\chi^1_\lambda(K_\nu)=1$,
$\chi^1_\lambda(K_\beta)=q^x$. Then the FIN condition implies that a
function $f\in H^0(\chi^1_{\lambda})$ is uniquely determined on elements of
the form $p_1(w)+p_2(w)F_\mu$ in ${\mathcal U}_q$, the elements
$p_1,p_2$ being polynomials in ${\mathcal A}_q$. Via the
correspondence (\ref{trans}), $H^0(\chi^1_{\lambda})$ is then
identified with the space ${\mathcal P}\otimes {\bm C}^2$, ${\mathcal
  P}$ being the usual space of polynomials on ${\bm C}^4$. It is now
easy to see that an element
\begin{equation}
  u^+_0=F^-_\mu F_\beta -(q+q^{-1})\cdot F_\beta E_\mu=w_2-q^{-1}\cdot
  w_1 F_\mu
\end{equation}
is an intertwiner into a space $H^0(\chi_{\delta})$ provided
$q^{2x-4}=1$. In this case, $\chi_{\delta}(K_\mu)=1,
\chi_{\delta}(K_\nu)=q$, and $\chi_{\delta}(K_\beta)=q^{x+1}$. In the
case where $x=2$ and $q$ is real there is a  unitary highest weight
representation in a subspace (the kernel of the operator $D^+$
below). If $x\neq 2$ we may still have an intertwiner provided that
$q=\E^{\frac{2\pi\cdot\I\cdot p}{2x-4}}$ for some integer $p$.

The target space thus consists of functions that are non-zero only on
expressions of the form $p_1(w)+p_2(w)F_\nu$. Now compute, for $f\in H^0(\chi^1_{\lambda})$,
\begin{eqnarray}
f_{u^+_0}(p_1(w)+p_2(w)F_\nu)&\kern-.7em=\kern-.7em&f([p_1(w)+p_2(w)F_\nu]
[w_2-q^{-1}w_1F_\mu])\\\nonumber
&\kern-.7em=\kern-.7em&f([p_1(w)w_2+p_2(w)w_4]-q^{-1}[p_1(w)w_1+ p_2(w)w_3] F_\mu).
\end{eqnarray}
Thus, the intertwiner $f\mapsto f_{u^+_0}$ may be expressed as a map
$D^+:{\mathcal P}\otimes {\bm C}^2\mapsto {\mathcal P}\otimes {\bm
  C}^2$ given by the matrix
\begin{eqnarray}\nonumber
  \label{eq:d+}
  D^+&=&\left(
    \begin{array}{cc}(w_2^R)^\dagger
      &(w_4^R)^\dagger\\-q^{-1}(w_1^R)^\dagger&-q^{-1}(w_3^R)^\dagger 
    \end{array} 
\right)\\&=&\left( 
    \begin{array}{cc}K_4\left[\frac{\partial}{\partial z_2}\right]_q
      &\left[\frac{\partial}{\partial
          z_4}\right]_q\\-q^{-1}K_2K_3K_4^2\left[\frac{\partial}{\partial 
  z_1}\right]_q&-q^{-1}K_4\left[\frac{\partial}{\partial z_3}\right]_q
    \end{array}
\right)\\\nonumber&+&\left( 
    \begin{array}{cc}0&0\\-q^{-1}z_4(1-q^{-2})K_4
      \square_q& 0
    \end{array}
\right).
\end{eqnarray}
Of course, there is an analogous operator $D^-$ obtainable by
interchanging $\mu$ and $\nu$;
\begin{eqnarray}\nonumber
  \label{eq:d-}
  D^-&=&\left(
    \begin{array}{cc}(w_3^R)^\dagger
      &(w_4^R)^\dagger\\-q^{-1}(w_1^R)^\dagger&-q^{-1}(w_2^R)^\dagger 
    \end{array} 
\right)\\&=&\left( 
    \begin{array}{cc}K_4\left[\frac{\partial}{\partial z_3}\right]_q
      &\left[\frac{\partial}{\partial
          z_4}\right]_q\\-q^{-1}K_2K_3K_4^2\left[\frac{\partial}{\partial 
  z_1}\right]_q&-q^{-1}K_4\left[\frac{\partial}{\partial z_2}\right]_q
    \end{array}
\right)\\\nonumber&+&\left( 
    \begin{array}{cc}0&0\\-q^{-1}z_4(1-q^{-2})K_4
      \square_q& 0
    \end{array}
\right).
\end{eqnarray} 
It is easy to see that $D^+D^-=D^-D^+=-q^{-1}\left( 
    \begin{array}{cc} \square_q&
      0\\0&\square_q 
    \end{array}
\right)$.

With the exception of the term $-q^{-1}z_4(1-q^{-2})K_4
\square_q$, the operators $D^\pm$ look like
any simple-minded generalization of the Dirac operator to the $q$
realm. However, there {\bf is} the extra term which represents the departure
from the ``quasi-commutative'' situation. We have in \cite{j7} seen
the first indications of an auxiliary bundle, which in fact can be
constructed using $\square_q$, in which
the new first order differential operators are the covariant
derivatives.  We leave this point to subsequent investigations.

\bibliographystyle{alpha} \bibliography{kernel}

\bigskip




\end{document}